\documentclass[a4paper,10pt,draft]{article}
\usepackage{babel,amsmath,amssymb,amsbsy,amsfonts,latexsym,amsthm}

\newtheorem{Def}{Definition}[section]
\newtheorem{teo}[Def]{Theorem}
\newtheorem{prop}[Def]{Proposition}
\newtheorem{lem}[Def]{Lemma}
\newtheorem{cor}[Def]{Corollary}
\newtheorem{oss}[Def]{Remark}
\newtheorem{example}[Def]{Example}

\newenvironment{Dim}[0]{{\bf Proof}} {}
\newenvironment{Dimo}[0]{{\bf{Proof of}}} 

\def\e{\varepsilon}

\def\o{\overline}
\def\u{\underline}
\def\w{\widehat}
\def\*{\star}

\def\t{\tilde}
\def\ra{\rangle}
\def\la{\langle}
\def\R{{\mathbb R}}

\def\N{{\mathbb N}}

\def\L1{\Lambda_B^1([0,T])}
\def\L2{\Lambda_B^2([0,T])}
\def\l{\lambda}
\def\L{\Lambda}

\def\M1{M^1_B([0,T])}
\def\M2{M^2_B([0,T])}

\def\I{{\cal I}}

\def\M{{\cal M}}
\def\N{{\cal N}}
\def\C{{\cal C}}

\def\I{{\cal I}}
\def\s{\sigma}
\def\p{\prime}
\def\u{\underline}
\def\o{\overline}
\def\t{\tilde}
\def\J{{\cal J}}
\def\P{{\cal P}}

\def\C{{\cal C}}
\def\S{{\cal S}}

\def\la{\langle}
\def\ra{\rangle}
\def\wx{\widehat{x}}
\def\wy{\widehat{y}}
\def\wt{\widehat{t}}
\def\var{\varepsilon}

\def\acc{\'}
\def\R{{\mathbb R}}

\newcommand{\cvd}{\begin{flushright}
	\rule[5pt]{5pt}{5pt}
	\end{flushright}}
\date{}
\begin{document}
\title{A Generalized Osgood Condition for Viscosity Solutions to Fully Nonlinear Parabolic Degenerate Equations}
\maketitle
\begin{center}
\author{Marco Papi\\
\vspace*{5pt}
\begin{scriptsize}
Istituto per le Applicazioni del Calcolo "M.Picone", V.le del Policlinico 137, I-00161 Roma (Italy) and\\ 
Dipartimento di Matematica di Roma "Tor Vergata", Via della Ricerca Scientifica, 00133, Roma (Italy),\\
e-mail: papi@iac.rm.cnr.it
\end{scriptsize}} 
\end{center}
\begin{center}
-----------------------------
\begin{quote}
\begin{small} ABSTRACT. - Using a generalized assumption of Osgood type,
we prove a new comparison result for viscosity sub and supersolutions
of fully nonlinear, possibly degenerate, parabolic equations.
Our result allows to deal with hamiltonian functions with a quadratic growth in the spatial
gradient, under special compatibility conditions with the diffusive terms.
It applies in particular to a financial differential model for pricing Mortgage-Backed Securities.
\end{small}
\end{quote}

-----------------------------
\begin{quote}
\begin{small}
\end{small}
\end{quote}
\end{center}
\section{Introduction}
In this paper we prove a new comparison result between viscosity sub and supersolution
for a nonlinear second-order parabolic, possibly strongly degenerate, equation of the 
following general form 
\begin{eqnarray}\label{GP}
\partial_t u+F(x,t,u,\nabla u,\nabla^2 u)=0,
\end{eqnarray}
in $\R^N\times[0,T)\;\; T>0$. The unknown $u$ will always be a real-valued function 
on $\R^N\times [0,T)$; $\partial_t u,\;\nabla u,\;\nabla^2 u$ denote respectively
the time derivative of $u$, the gradient of $u$ and the Hessian of $u$ in the space variables. 
$F$ is
a real-valued function defined on 
$\R^N\times[0,T)\times(a,b)\times\R^N\times\S^N$, where $\S^N$ is the space of 
$N\times N$ symmetric matrices 
endowed with the usual ordering, while $(a,b)$ is an open, possibly unbounded interval.\\
Our comparison principle is based on a special compatibility condition between the diffusive terms and the
quadratic dependence of the Hamiltonian on the spatial gradient, see formula (\ref{CP7}) in Theorem \ref{CP}.
This condition is somehow similar to the celebrated ``null condition", used to prove global existence for hyperbolic
equations with a quadratic dependence on the gradient, see for instance \cite{Hormander}, \cite{KK} and references therein. Notice
however that with our comparison result, the  dependence of the Hamiltonian on the unknown $u$ is allowed.
Moreover, it includes not only the  quadratic growth
with respect to the spatial gradient, but also the absence of monotonicity with respect to $u$ and a lack
of regularity in the Hamiltonian.\\
Actually one important motivation which has brought us to state this principle, comes from the following quasilinear equation:
\begin{eqnarray}\label{DM1}
\!\partial_t U\!-\!\frac{1}{2}tr(\s\s^{\top}\nabla^2 U)\!-\!\la \mu,\nabla U \ra\!+\!\rho
\frac{|\s^{\top}\nabla U|^2}{U+h+\xi(t)}
+r(U\!+\!h)\!-\!\tau h=0,
\end{eqnarray}
in $\R^N\times[0,T),\;\rho>0,\;\tau,\;T>0$, where $tr$, $\la,\ra$, $|\cdot|$ denote the 
trace of a square matrix, the Euclidean norm and inner product, respectively, and with 
a con\-ti-nuous initial datum $U_0$. 
Moreover 
$\mu:\R^N\times[0,T)\rightarrow \R^N$, $\s:[0,T)\rightarrow \M_{N\times d}(\R)$, 
$\xi,\; r:[0,T)\rightarrow [0,\infty)$,
$h:\R^N\times[0,T)\rightarrow \R$ are continuous, where $\M_{N\times d}(\R)$ 
denotes the space of real $N\times d$
matrices, with $N\geq d$.\\
Equation (\ref{DM1}) has been proposed in \cite{37}, following the probabilistic financial model by X. Gabaix and O. Vigneron \cite{19}, as a differential model for pricing some widely used American financial instruments, the Mortgage-Backed 
Securities ($MBS$). Although we shall not study here the financial issues coming from 
model (\ref{DM1}), in section \ref{MOTIVE} 
we quickly analyse the  particular structure 
of equation (\ref{DM1}).\\ 
Quasilinear parabolic and elliptic 
equations have been extensively studied in the literature in the 
uniformly or strictly 
elliptic case. We can quote some classical books for a wide  presentation
of the results 
which are known in this direction, \cite{TR}, \cite{LAD}. For real applications
we need to consider $N>d$, but the case of strong de\-ge\-ne\-ra\-te equations, like 
(\ref{DM1}), is far less classical, since we cannot expect the global existence of strong 
solutions. It turns out that viscosity 
solutions can be used to treat this kind of problem, see \cite{1}. For stationary degenerate elliptic, but also 
parabolic equations, when
we want to prove a comparison result by using viscosity solution's method, the 
coercivity of the equation
in the variable $u$ is a standard assumption, namely for any $R>0$, there 
exists $\gamma_R$, which
is positive in the stationary case and possibly negative in the parabolic case,
such that for any $-R\leq v \leq u\leq R$, $p\in\R^N$ and $X$, symmetric 
matrix with real coefficients, it holds
\begin{eqnarray}\label{8}
F(u,p,X)-F(v,p,X)\geq \gamma_R (u-v).
\end{eqnarray}
This assumption can be more complicated to be verified if there is dependence by the  spatial variable and the
domain of $F$ with respect to $u$ is not all of $\R$. Condition (\ref{8}) 
is then used to state comparison results
also in unbounded domains such as in \cite{3}. There are many ways to relax 
condition (\ref{8}); for instance 
in \cite{22} a Osgood condition on $F$ is used (see Definition \ref{OSGOOD}), but also 
this condition
is not enough for equation (\ref{DM1}), because the dependence on the term $(U+h+\xi)^{-1}$. Assumptions like (\ref{8}) can be removed, if it is possible to prove that the 
Hopf's Lemma
holds in the strong form. The problem with this approach is that, for quasilinear 
equations, the Hopf's Lemma 
can be obtained only for classical solutions. In a recent work of G. Barles and J. Busca 
\cite{Barles}, a
weak version of the Hopf's Lemma is used to show the comparison in bounded domain for a 
stationary fully nonlinear equations, but
because of the de\-ge\-ne\-ra\-cy of equation (\ref{DM1}) it is not possible to follow here the same 
technique. 
Our proof of the comparison principle for (\ref{GP}), proved in section \ref{COMPARISON}, follows a different idea, based on the structural condition (\ref{CP7})
on the function $F$.\\
Finally, let us notice that a similar problem has been recently considered by P.L. Lions and P.E. Souganidis in \cite{LS}. 
Starting from a fully nonlinear stochastic partial differential equation, they reduce their problem to a deterministic second order fully nonlinear pde and introduce some conditions for the 
original Hamiltonian to obtain the comparison principle.\\
In section \ref{COMPARISON}, we give some examples to show that the conditions proposed in that paper are strictly contained in our approach.
Moreover the present results on equation (\ref{DM1}) do not follow by a simple combination of the ideas
in \cite{LS}, see Remark \ref{MBSmodel}.\\
In a work in preparation \cite{spde}, following the general theory of stochastic partial differential
equations of \cite{LS17}, \cite{LS18}, \cite{LS} we use the same kind of conditions on the Hamiltonian, here proposed in our Theorem \ref{CP},
to proving an existence result for stochastic pde, which depend also on $(x,t,u)$.

\section{The MBS Model}\label{MOTIVE}
In the last years, the theory about financial markets, the mathematical frameworks
for
 modelling them, and arbitrage theory have reached an high de\-ve\-lop\-ment and have taken
a prominent position in the mathematical literature (we refer the interested reader
to some reference books about these arguments like \cite{17}, \cite{18}). However some  
financial derivatives, like MBSs, still need of specific mo\-dels which are conformed on their 
peculiarities. Morover when one introduces the incompleteness in the market, 
then the choice of a new measure to e\-va\-lua\-te the derivatives
is necessary. Via the Girsanov Theorem, \cite{9}, this choice corresponds to assign
a \textit{market price of risk}.
Essentially when the market is not complete, the market price of risk is not unique and 
it must be chosen by statistical
analysis; this approach is also followed by X.Gabaix and O.Vigneron in \cite{19}, 
but the mathematical
toolbox that would required  to rigorously derive their results is still largely to be 
developed. In \cite{37}, we used arbitrage arguments and the form of the market price 
of risk proposed in \cite{19} to develop our differential model (\ref{DM1}).
The appeal to differential equations in financial mo\-del\-ling 
has become a standard approach, and in many cases it represents the best way for valuing
derivatives. In some instances, the expectations 
that result from applying the arbitrage pricing principle
can be characterized as solutions to partial differential equations. Therefore 
the study of solutions of the pricing equations of more complex financial 
instruments and their numerical solutions have become 
important techniques available to practitioners of modern quantitative  
finance. However up to now the $PDE$ approach for valuation of MBSs with the additional
specification of a nonlinear form for the market price of risk, was not followed. Although 
the equation (\ref{DM1}) represents a reduced version of a more general si\-tua\-tion which includes 
a model of the issuance of future securities, 
the model outlines a new manner for treating MBSs derivatives. Not only the possibility to 
represents the value of an MBS
as the solution of a $PDE$ is important, but also the existence 
and uniqueness problem, related to that equation holds
an essential role: actually it is well known that an arbitrary specification of the market 
price of risk, may lead to arbitrage opportunities, hence the existence of a solution for 
(\ref{DM1}) is the proof that exists a new $risk-neutral$ market measure through evaluating these securities. This equation describes a first model in this direction, actually 
it is considered a deterministic interest rate, while in a future work, at the moment in 
preparation, \cite{39}, we study a model which takes into account different stochastic models for interest rate movements. However the model receives a reasonable 
empirical support.\\
Mortgage-Backed derivatives
are the products of a securitization of pools of mortgages. Every mortgage-holder in a pool
holds the right to prepay her debt at every time between $0$ and the maturity $T$. This 
American-style $option$ determines the stochasticity in 
the valuation of the price of a Mortgage-Backed
from that pool. Many factors affect prepayments and this fact creates a 
remarkable complexity
for generating forecasts.\\
The arbitrage pricing principle applies to the financial instruments whose cash
flows are related to the values of economic factors, such as the interest rates ($r=r(t)$). It implies that the
price ($V_t$) of such derivative assets can be expressed as a conditional expectation 
with respect to a particular measure ($Q$), over the probability space of the 
underlying factors that affect its value. Under
suitable conditions, computing this expectation it
reduces to the solution of a partial differential equation:
in fact $V_t=U(X_t,T-t)+h(X_t,T-t)$, where $h$ contains information about 
prepayments and $X_t$ describes
economic factors. In \cite{19} the measure $Q$ depends in a 
nonlinear way by the price $V$ and its vo\-la\-ti\-li\-ty.
This aspect produces the nonlinear quadratic term in 
the equation (\ref{DM1}). Another characteristic of problem (\ref{DM1})
is the strongly degeneracy of the equation. In fact, also this aspect comes 
from the dependence of the payoff of the
security on the trajectory followed by one or more of the underlying markovian processes.
We illustrate this point with a simple consideration. One possible 
index ($b_t$) of the incentive to prepay a mortgage is the amount by which
a particular interest rate $\o{r}_s$ is below some given level $\o{r}$. The factor
is represented by $ b_t =\int_0^t (\o{r}-\o{r}_s)^+ ds $, and this form gives 
the absence of diffusion
in the $b$ direction. 

\section{Comparison Principle and Existence}\label{COMPARISON}
We begin with some notation and recall about viscosity solutions.\\
Denote $Lip(\R^N)$, $\C^{2,1}(\R^N\times [0,T))$, $\C^{k}(I)$, respectively,
the space of the globally Lipschitz functions over $\R^N$, the space
of all functions which have two continuous derivatives in the space variable and one continuous in
the time over $\R^N\times[0,T)$, the space of all functions with $k=0,\;1,\;2$ continuous derivatives
over the interval $I$.
\begin{Def}\label{P2+}
Given an upper semicontinuos function $u:\R^{N}\times [0,T)\rightarrow \R$, the parabolic super 2-jet  
of $u$ at the point $(x,t)$ is,
\begin{eqnarray*}
  &&\!\!\!{\cal P}^{2,+}u(x,t)\!\!:=\!\!\Big\{(\partial\varphi(x,t),\nabla\varphi(x,t),\nabla^{2}\varphi(x,t))
     :\;\varphi\in{\cal C}^{2,1}(\R^{N}\times[0,T)),
      u\!\!-\!\!\varphi \!\!\!\\ &&\;has\;a\;local\;maximum\;at\;(x,t)\Big\}
\!\!=\!\!\Big\{(b,p,X)\in\R\times\R^N\times \S^N:\;
     u(y,s)\leq \\ && u(x,t)\!+\!b(t-s)\!+\!\la p, y-x\ra\!+\!\frac{1}{2}\la X(y-x),(y-x)\ra\!+\!o(|y-x|^2+|s-t|),\;\;
	\\ && as\;\;\!(y,s)\rightarrow (x,t)\Big\}.
\end{eqnarray*}
\end{Def}
Similarly, if $u$ is a lower semicontinuous function, we define the lower 2-jet as $\P^{2,-}u=-\P^{2,+}(-u)$.

\begin{oss}
We shall say that a function $\varphi$, as in the previous definition, is a $test$ $function$ for $ {\cal P}^{2,+}u$ at $(x,t)$.\\
In the definition (\ref{P2+}) it is possible to replace local by global, or local strict, or global strict. 
In the case of a global strict maximum
with $(u-\varphi)(x,t)=0$ we will say that $\varphi$ is a \textit{good test function} 
for $ {\cal P}^{2,+}u$ at $(x,t)$.\\
\end{oss}


Denote with $u_{\*}$, $u^{\*}$ respectively the lower and upper semicontinuous envelope
of $u$. Moreover we consider as domain of $F$ with respect to $u$, in the equation (\ref{GP}),
a (possibly bounded) open interval $(a,b)$, and as initial datum $u_0\in \C^0((a,b))$.
\begin{Def}\label{sottosoluzione}
A function $u:\R^N\times[0,T)\rightarrow \R$, is called a viscosity sub (resp. super) solution 
of (\ref{GP}) if $a<u\leq u^{\*}<b$ (resp. $a<u_{\*}\leq u<b$), and 
\begin{description}
     \item[1.]for every $(x,t)\in \R^{N}\times [0,T)$ and a test function $\varphi$ for
    ${\cal P}^{2,+}u^{\*}$ (resp. ${\cal P}^{2,-}u^{\*}$) at $(x,t)$,  
     \[ \partial_{t}\varphi(x,t)+F(x,t,u^{\*}(x,t),\nabla \varphi(x,t),\nabla^{2}
    \varphi(x,t))\leq 0\]
    \[ (resp.\;\;\; \partial_{t}\varphi(x,t)+F(x,t,u_{\*}(x,t),\nabla \varphi(x,t),\nabla^{2}
    \varphi(x,t))\geq 0
     ),\]
     \item[2.]$u(x,0)\leq u_{0}(x)$, (resp. $u(x,0)\geq u_{0}(x)$) for $x\in\R^N$.
\end{description}
\end{Def}

\begin{Def}
A function $u:\R^N\times[0,T)\rightarrow \R$ is called a viscosity solution of (\ref{GP}) 
if it is, at the same time, a viscosity sub- and a super- solution of (\ref{GP}).
\end{Def}


In this paper we call a continuous function $\nu:[0,\infty)\rightarrow [0,\infty)$ a modulus, if
$\nu(0)=0$ and it is nondecreasing; given $p\in\R^N$, we denote $p\otimes p$ the $N\times N$ matrix whose entries are $p_i p_j$ for every $i,\;j=1,\ldots,N$; if $X\in\S^N$
then $\|X\|$ denotes the operator norm of $X$ as a self-adjoint operator on $\R^N$; if $f$ is a real-valued bounded function in a domain
$D$, then we denote with $\|f\|_{\infty}$ the supremum of $f$ over $D$.
Moreover we need of the following definition in which we recall a 
property already used in \cite{22}.

\pagebreak

\begin{Def}\label{DOSGOOD}
We say that a continuous function $\Gamma:[0,l]\rightarrow \R$, $l>0$ satisfies the {\em Osgood Condition}, if
the following conditions hold:
\begin{description}
\item[($i$)] $\Gamma$ is an increasing function and $\Gamma(0)=0$;
\item[($ii$)] 
\[\int_0^l \frac{dr}{\Gamma(r)}=+\infty.\]
\end{description}
\end{Def}

We are now in position to state our main comparison Theorem.

\begin{teo}\label{CP}
Consider the following differential equation:
\begin{eqnarray}\label{CP1}
\partial_t u+F(x,t,u,\nabla u, \nabla^2 u)=0,\;\;\;\;\;(x,t)\in \R^{N}\times(0,T)
\end{eqnarray}
and assume that,
	\begin{eqnarray}\label{CP2}
		 F\in\C(\R^{N}\times[0,T)\times(a-\var_0,b+\var_0)\times\R^N\times\S^{N}),\;\;\;\;\var_0>0;
	\end{eqnarray}
- $F$ is degenerate elliptic, i.e.,
	\begin{eqnarray}\label{CP3}
		F(x,t,u,p,X+Y)\leq F(x,t,u,p,X),
	\end{eqnarray}
	for every $(x,t,u,p)\in\R^{N}\times[0,T)\times(a-\var_0,b+\var_0)\times\R^N, X,\;Y\in\S^N , Y\geq 0$.\\
- For every $R>0$ there is a modulus $\nu_{1,R}$ such that 
	\begin{eqnarray}\label{CP4}
		|F(x,t,u,p,X)-F(x,t,u,q,X)|\leq \nu_{1,R}(|p-q|)
	\end{eqnarray}
	for all $(x,t,u,p,X)\in\R^{N}\times[0,T)\times[a,b]\times\R^N\times\S^{N}$, with $|p|,|q|,\|X\|\leq R$.\\
- Suppose that,
	\begin{equation}\label{CP5}
		-\varepsilon_1 I_{2N}\leq\left(\begin{array}{cc}
					X & O\\
					O & Y
			\end{array}	\right)\leq \varepsilon_2 
                                          \left(\begin{array}{cc}
					I_N & -I_N\\
				       -I_N &  I_N
			\end{array}	\right)  +\varepsilon_3 I_{2N}
	\end{equation}
	with $X,\;Y\in\S^N$, $\varepsilon_1,\e_2,\e_3\geq 0$. Let $R$ be taken so that $R\geq max(\e_1,2\e_2+\e_3)+2\e_3$. 
	Then it holds:
	\begin{eqnarray}\label{CP6}
		\lefteqn{F(x,t,u,p,X+Z)-F(y,t,u,p,-Y+Z) \geq \nonumber }\\
  			&&-\nu_{2}\big(|x-y|(|p|+1)+\e_2 |x-y|^2\big)
						-\nu_{2,R}(2\e_3)
	\end{eqnarray}  
	for every $(x,t,u)\in\R^N\times[0,T]\times[a,b]$, $|p|,\;\|Z\|\leq R$, and with some moduli $\nu_2,\;\nu_{2,R}$ indipendent of the other variables and $\e_1,\;\e_2,\;\e_3$; where
	$\nu_{2,R}$ possibly dipendent of $R$.\\
- There are functions $\Gamma$ which satisfies the {\em Osgood Condition} \ref{DOSGOOD} over \\ $[0,\sqrt{\frac{\L_0}{\l_0}}(b-a)]$, and $z\in\C^1((a-\e_0,b+\e_0);(0,+\infty))$, with $z([a,b])\subset [\l_0,\L_0]$, $\L_0>\l_0>0$,
	such that for every $R>0$ there is a modulus $\w{\nu}_R$, such that,
 	\begin{eqnarray}\label{CP7}
		&&\frac{1}{\l}F(x,t,u,\l q, \l X +\kappa q\otimes q)-\frac{1}{\w\l}F(x,t,v,\w\l q, \w\l X+\w\kappa q\otimes q)\geq\nonumber\\
		&&\geq-\Gamma(u-v)-\w{\nu}_R\big( (|\l^2-z(u)|+|{\w\l}^2-z(v)|)(1+|q|+\|X\|)\big)
	\end{eqnarray}
	for $(x,t)\in\R^N\times[0,T)$, $|q|,\;\|X\|\leq R$, $a \leq v\leq u\leq b $, $\inf_{[a,b]} z^{\frac{1}{2}}\leq \l,\;\w\l\leq\sup_{[a,b]}z^{\frac{1}{2}}$, $2\kappa\leq z^{\p}(u)$, $2\w\kappa\geq z^{\p}(v)$.\\
If $\u{u},\;\o{u}:\R^N\times[0,T)\rightarrow [a,b]$, are respectively viscosity subsolution and supersolution of 
the equation (\ref{CP1}) and moreover there is a 
modulus $\omega$ such that
\begin{eqnarray}\label{CP8}
\u{u}^{\star}(x,0)-\o{u}_{\star}(y,0)\leq \omega(|x-y|)
\end{eqnarray} 
for all $x,y\in\R^N$, then 
\begin{eqnarray}\label{CP16}
	\u{u}^{\star}\leq \o{u}_{\star},
\end{eqnarray}
over $\R^N\times[0,T)$.
\end{teo}

\begin{oss}
As regards the bounds of the solution, by financial purposes
we look for bounded solution, so we limit us to state the comparison
only for bounded sub/supersolution which take values only in
subsets of the domain of $F$.\\
Conditions (\ref{CP3})-(\ref{CP6}) were already used in \cite{3}, but
assumption (\ref{CP7}), is a relaxation of ``monotonicity" assumption
with respect $u$, about $F$. \end{oss}
We shall prove Theorem \ref{CP} in some steps. We start by describing
some technical results which will be useful in the proof of Theorem \ref{CP}.
For the proof of these we refer the reader to the works which contain them.


For parabolic problems the monotonicity assumption with respect to $u$ can be relaxed
requiring the Osgood type condition on $F$. Actually the proof of Theorem is based
on the following Lemma. 

\begin{lem}(\cite{22})\label{OSGOOD}
Let $\Gamma$ be an Osgood type function over $[0,l]$, in the sense of Definition \ref{DOSGOOD}. 
Let $f$ be an upper semicontinuous function over $[0,T)$, valued in $[0,l]$. Assume that $f$ satisfies in the viscosity sense the following problem,
\begin{equation}\label{l}
	\left\{\begin{array}{ll}
		f^{\prime}(t)\leq \Gamma(f(t)),  & \mbox{$t \in (0,T)$}\\\\
		\min(f^{\prime}(0)-\Gamma(f(0)),f(0))=0.
               \end{array}
	\right.
\end{equation}
Then $f\equiv 0$.
\end{lem}

The last Lemma is an important Technical Lemma, and is Proposition $4$ in \cite{22}.

\begin{Def}
The super 1-jet of an upper semicontinuos function \\$u:[0,T)\rightarrow\R$ at the point $t_0$ is
\begin{eqnarray*}
    {\cal J}^{1,+}u(t_0):= \Big\{\varphi^{\p}(t_0):\;\varphi\in {\cal C}^{1}([0,T)),\;\;u-\varphi\;has\;a\;
	 local\;maximum \;at\;t_0\Big\}
\end{eqnarray*}
\end{Def}

\begin{lem}\label{LEMMATECNICO}(\cite{22}, Proposition 4)
Let $\underline{v}$, $\overline{v}$ be respectively an upper/lower semicontinuous function, such that 
\begin{eqnarray*}
\	\sup\Big\{(\underline{v}(x,t)-\overline{v}(y,t))^{+}\;:\;|x-y|\leq 1\Big\}\leq K
\end{eqnarray*}
Set the function, 
\begin{eqnarray*}
	\vartheta(t)=\lim_{r\rightarrow 0}\sup\Big\{(\underline{v}(x,t)-\overline{v}(y,t))^{+}
	\;:\;|x-y|\leq r\Big\}.
\end{eqnarray*}
Denote with $\vartheta^{\star}$ the upper semicontinuous envelope of $\vartheta$, and let $\varphi$ 
be a good test function for ${\cal J}^{1,+}\vartheta^{\star}(t_0)$.\\
Consider the function (which depends on two positive parameters $\alpha,\delta$)
\begin{equation}\label{Phi}
	\Phi_{\alpha,\delta}(x,y,t)=\underline{v}(x,t)-\overline{v}(y,t)-\frac{\alpha}{2}
	|x-y|^{2}-\delta|x|^{2}-\varphi(t).
\end{equation}
Set $\Delta=\{(x,y,t)\;:\;|x-y|\leq 1,\;t\in[0,\frac{t_0+T}{2}]\}$, and 
$\overline{\Phi}_{\alpha,\delta}=\sup_{\Delta}\Phi_{\alpha,\delta}$.\\
Then 
\begin{description}
	\item[$(i)$] for every $\alpha,\delta$ exists $(x_{\alpha,\delta},y_{\alpha,\delta},
	t_{\alpha,\delta})\in \Delta$ which is a maximum point of $\Phi_{\alpha,\delta}$; \\
	\item[$(ii)$]$\lim_{\delta\rightarrow 0}\overline{\Phi}_{\alpha,\delta}=
	\sup_{\Delta}\{\Phi_{\alpha,0}\}\equiv\overline{\Phi}_{\alpha}$;\\
	\item[$(iii)$]$\lim_{\alpha\rightarrow +\infty}\overline{\Phi}_{\alpha}=
	\vartheta^{\star}(t_0)-\varphi(t_0)=0$;\\


	\item[$(iv)$]$\lim_{\alpha\rightarrow +\infty}\lim_{\delta\rightarrow 0}
	(\alpha|x_{\alpha,\delta}-y_{\alpha,\delta}|^{2}+\delta|x_{\alpha,\delta}|^{2})=
	0$.
\end{description}
Moreover, for a subsequence we can obtain,
\begin{description}
	\item[$(v)$]$\lim_{\alpha\rightarrow +\infty}\lim_{\delta\rightarrow 0}
	t_{\alpha,\delta}=t_0$;\\
	\item[$(vi)$]$\lim_{\alpha\rightarrow +\infty}\lim_{\delta\rightarrow 0}
	\Big(\underline{v}(x_{\alpha,\delta},t_{\alpha,\delta})-
	\overline{v}(y_{\alpha,\delta},t_{\alpha,\delta})\Big)=\vartheta^{\star}(t_0)=
	\varphi(t_0)$.
\end{description}
\end{lem}

Before to give the proof of Theorem \ref{CP}, we recall the well known pro\-per\-ty about the conservation of the notion of
viscosity solution with respect to a change of variable of the unknown, with respect to an increasing smooth transformation.\\

\begin{Dimo} {\bf Theorem \ref{CP}.} Consider the application, $\Psi:(a-\frac{\textstyle{\e_0}}{\textstyle{2}},b+\frac{\textstyle{\e_0}}{\textstyle{2}})
\rightarrow \R$, defined by
\begin{eqnarray*} \Psi(u)=\int_{a-\frac{\textstyle{\e_0}}{\textstyle{2}}}^{u}
\frac{1}{\sqrt{z(\tau)}}d\tau,\;\;\;\;u\in\Big(a-\frac{\textstyle{\e_0}}{\textstyle{2}},b+\frac{\textstyle{\e_0}}
{\textstyle{2}}\Big) 
\end{eqnarray*} since that
$z>0$, $\Psi$ is well defined and $\Psi^{\p}>0$, so $\Psi$ has a $\C^2$
inverse which we denote with $\I$ defined in $(\Psi(a-\frac{\textstyle{\e_0}}{\textstyle{2}}),\Psi(b+\frac{\textstyle{\e_0}}
{\textstyle{2}}))$. 
If we consider the functions
\begin{eqnarray*} \u{v}=\Psi\circ \u{u},\;\;\;\;\;\;\;,\o{v}=\Psi\circ
\o{u}, \end{eqnarray*}
then by the increasing property of $\Psi$, $\u{v},\;\o{v}$ are respectively subsolution and supersolution of the
following equation, 
\begin{eqnarray}\label{CP9} 
\partial_t v+\t{F}(x,t,v,\nabla v,\nabla^2 v)=0,\;\;\; (x,t)\in\R^N\times(0,T)
\end{eqnarray} \begin{eqnarray}\label{CP10}
\t{F}(x,t,v,p,X)=\frac{1}{\I^{\p}(v)}F(x,t,\I(v),\I^{\p}(v)p,
\I^{\p}(v)X+\I^{\p\p}(v)p\otimes p), \end{eqnarray} for every
$(x,t,v,p,X)\in\R^m\times[0,T)\times (\Psi(a-\frac{\textstyle{\e_0}}{\textstyle{2}}),\Psi(b+\frac{\textstyle{\e_0}}
{\textstyle{2}}))\times\R^N\times\S^N$
and, by (\ref{CP8}), \begin{eqnarray}\label{CP11}
\u{v}^{\star}(x,0)-\o{v}_{\star}(y,0)\leq \frac{1}{\sqrt
\l_0}\omega(|x-y|),\;\;\;\;x,y\in\R^N. \end{eqnarray} we will prove the
comparison between $\u{v}^{\star},\;\o{v}_{\star}$. Set
\[\vartheta(t):=\lim_{r\rightarrow
0^+}\sup\Big\{(\u{v}^{\star}(x,t)-\o{v}_{\star}(y,t))^+:\; |x-y|\leq
r\Big\}\] and denote $\vartheta^{\star}\leq \frac{b-a}{\sqrt{\l_0}}$ its upper-semicontinuous
envelope. Our aim is to show that $\vartheta^{\star}\equiv 0$.\\ We will
obtain this assertion, proving that \[\left\{\begin{array}{ll}
(\vartheta^{\star})^{\prime}(t)\leq \Gamma_0
(\vartheta^{\*}(t))&\mbox{$t\in (0,T)$}\\\\
\min((\vartheta^{\star})^{\prime}(0)-\Gamma_0 (\vartheta^{\*}(0)),\:
\vartheta^{\*}(0))=0 \end{array}\right.\] where, $\Gamma_0(\theta):=\Gamma(\sqrt{\L_0}\theta)$ and  
$\Gamma$ is given by(\ref{CP7}), it holds in a viscosity sense, and then using Lemma
\ref{OSGOOD}, in fact it easy to see that $\Gamma_0$ is an Osgood type function over $[0,\frac{b-a}{\sqrt{\l_0}}]$.\\ Let $t_0 \in [0,T)$: \begin{description} \item[-] if
$\vartheta^{\star}(t_0)=0$ and $t_0=0$, we have nothing to say .
\item[-] if $\vartheta^{\star}(t_0)=0$ and $t_0>0$, $\vartheta^{\star}$
has an interior minimum at $t_0$ and this is the same for any test
function $\varphi$ for $\J^{+}\vartheta^{\star}(t_0)$. Since that $\varphi$
is regular, follow that,
\[\varphi^{\prime}(t_0)=0=\Gamma_0(0)=\Gamma_0(\vartheta^{\star}(t_0))\]
\item[-]if $\vartheta^{\star}(t_0)>0$, let $\varphi$ be a good test
function for $\J^{1,+}\vartheta^{\star}(t_0)$. By the boundness of $\u{v}$
and $\o{v}$, we can apply Lemma \ref{LEMMATECNICO}. Let
$\Phi_{\alpha,\delta}$ be the function defined in $(\ref{Phi})$; for
fixed $\alpha,\delta$, this function obtain a maximum in a point which
we denote with $(z_1,z_2,s)$ (leaving out the dipendence from the
parameters, in order to simplify the notations).\\ Moreover, since that
$\vartheta^{\star}(t_0)>0$ and by the conditions $(iv)$, $(v)$, $(vi)$
of Lemma \ref{LEMMATECNICO} and (\ref{CP11}), definitively, for large
$\alpha$ and small $\delta$, we can assume
$\u{v}^{\*}(z_1,s)>\o{v}_{\*}(z_2,s)$, $\delta z_1\leq \sqrt{\delta}$, $s>0$ . 
\end{description}

In order to obtain information about the ``derivatives" of $\vartheta^{\star}$, and noting that $s \in (0,T)$, we can apply 
the classical maximum principle for semicontinuous functions, which is Theorem $6$ in \cite{2},
\[\begin{array}{lr}
	u_1(x,t)=\u{v}^{\*}(x,t)-\delta|x|^2,\;\; & u_2(y,t)=-\o{v}_{\star}(y,t),
\end{array}\]
\[w(x,y,t)=u_1(x,t)+u_2(y,t),\]
for every $(x,y,t)\in \R^N\times\R^N\times (0,T)$. Set
\[b=\varphi^{\prime}(s),\;p=\left(\begin{array}{c}
		\alpha(z_1-z_2)\\
		-\alpha(z_1-z_2)
	\end{array}\right) ,\;A=\alpha\left(\begin{array}{cc}
				I_N & -I_N\\
				-I_N & I_N
				\end{array}\right)\]
where $I_N$ is the $N\times N$ identity matrix. Then 
$(b,p,A)\in \P^{2,+}w(z_1,z_2,s)$.\\
There exists $(b_1,X_1),(b_2,X_2)\in \R\times\S^N$ 
	such that
	\[(b_1,\alpha(z_1-z_2),X_1)\in \o{\P}^{2,+}(\u{v}^{\*}-\delta|\cdot|^2)(z_1,s)\]
	and
	\[(b_2,-\alpha(z_1-z_2),X_2)\in \o{\P}^{2,+}(-\o{v}_{\*})(z_2,s),\]


        so,	
	\[(b_1,\alpha(z_1-z_2)+2\delta z_1,X_1+2\delta I_N)\in \o{\P}^{2,+}\u{v}^{\*}(z_1,s)\]
	\[(-b_2,\alpha (z_1-z_2),-X_2)\in \o{\P}^{2,-}\o{v}_{\*}(z_2,s).\]
	Moreover, since that $A^{2}=2 \alpha A$, and $\|A\|=2\alpha$, the following relations hold
\begin{equation}\label{CP12}
-3 \alpha I_{2N} \leq \left(\begin{array}{cc}
				X_1 & 0\\
				0 & X_2
				\end{array}\right)\leq3\alpha
                                \left(\begin{array}{cc}
				I_N & -I_N\\
				-I_N & I_N
				\end{array}\right)
\end{equation}
\begin{equation}\label{CP13}
	b_1+b_2=\varphi^{\prime}(s)
\end{equation}
where $I_{2N}$ denotes the $2N\times 2N$ identity matrix.
We observe that, for continuity reasons, the equation (\ref{CP9}) is also satisfied, by $\u{v},\o{v}$,
over the closure of the parabolic-semijets. So we have,

\begin{eqnarray*}
	b_1\leq&-&\t{F}(z_1,s,\u{v}^{\*}(z_1,s),\alpha(z_1-z_2)+2\delta z_1,X_1+2 \delta I_m),\\\\
	b_2\leq&+&\t{F}(z_2,s,\o{v}_{\*}(z_2,s),\alpha(z_1-z_2),-X_1).
\end{eqnarray*}
 
Set $p:=\alpha(z_1-z_2)$, then these two relations and (\ref{CP13}) 
imply,
\begin{eqnarray}\label{CP15}
	\varphi^{\prime}(s) &\leq& \t{F}(z_2,s,\o{v}_{\*}(z_2,s),p,-X_2)+\nonumber\\
	&&-\t{F}(z_1,s,\u{v}^{\*}(z_1,s),p+2\delta z_1,X_1+2 \delta I_N).
\end{eqnarray}
\begin{eqnarray*}
&&\!\! \t{F}(z_2,s,\o{v}_{\*}(z_2,s),p,-X_2)
\!-\!\t{F}(z_1,s,\u{v}^{\*}(z_1,s),p+2\delta z_1,X_1+2 \delta I_N)=\\\\
&& \!\Big[\! \t{F}(z_2,s,\o{v}_{\*}(z_2,s),p,-X_2)\!\!-\!\frac{1}{\I^{\p}(\o{v}_{\*}(z_2,s))}\!
 F(z_2,s,\!\I(\o{v}_{\*}(z_2,s)),\!\I^{\p}(\o{v}_{\*}(z_2,s))\times\\\\
&& (p+2\delta z_1),
-\I^{\p}(\o{v}_{\*}(z_2,s))X_2+\I^{\p\p}(\o{v}_{\*}(z_2,s))p\otimes p)\Big]\!
\!+\!\frac{1}{\I^{\p}(\o{v}_{\*}(z_2,s))}\times \\\\
&&\Big[F(z_2,s,\I(\o{v}_{\*}(z_2,s)),\I^{\p}(\o{v}_{\*}(z_2,s))(p+2\delta z_1),
-\I^{\p}(\o{v}_{\*}(z_2,s))\times \\\\
&& X_2+\I^{\p\p}(\o{v}_{\*}(z_2,s))p\otimes p)\\\\
&&-F(z_1,s,\I(\o{v}_{\*}(z_2,s)),\I^{\p}(\o{v}_{\*}(z_2,s))(p+2\delta z_1),
\I^{\p}(\o{v}_{\*}(z_2,s))(X_1+2\delta I_N)+\\\\
&&\I^{\p\p}(\o{v}_{\*}(z_2,s))(p+2\delta z_1)\otimes(p+2\delta z_1))\Big]\\\\
&&+\Big[\frac{1}{\I^{\p}(\o{v}_{\*}(z_2,s))}F(z_1,s,\I(\o{v}_{\*}(z_2,s)),\I^{\p}(\o{v}_{\*}(z_2,s))(p+2\delta z_1),
\I^{\p}(\o{v}_{\*}(z_2,s))\times\\\\
&&(X_1+2\delta I_N)+
\I^{\p\p}(\o{v}_{\*}(z_2,s))(p+2\delta z_1)\otimes(p+2\delta z_1))-
\frac{1}{\I^{\p}(\u{v}^{\*}(z_1,s))}\times \\\\ && F(z_1,s,\I(\u{v}^{\*}(z_1,s)),\I^{\p}(\u{v}^{\*}(z_1,s))(p+2\delta z_1),
\I^{\p}(\u{v}^{\*}(z_1,s))(X_1+2\delta I_N)+\\\\
&&\I^{\p\p}(\u{v}^{\*}(z_1,s))(p+2\delta z_1)\otimes(p+2\delta z_1))\Big]
\end{eqnarray*}
Now we estimate the single terms in the brackets $[\cdot ]$ in the right hand side of this equality.
For sufficiently large $\alpha$, and $\delta<1$, set $R_{\alpha}:=(\alpha+2)^2\max(\sqrt{\L_0},\frac{\textstyle{\|z^{\p}\|_{\infty}}}{\textstyle{2}})$\\ $+3\alpha\sqrt{\L_0}$ in (\ref{CP4}), then by (\ref{CP10}), the inequalities (\ref{CP12}),
Lemma \ref{LEMMATECNICO} and the assumptions on the function $z$, we infer,
\begin{eqnarray*}
	&& \t{F}(z_2,s,\o{v}_{\*}(z_2,s),p,-X_2)-\frac{1}{\I^{\p}(\o{v}_{\*}(z_2,s))}F(z_2,s,\I(\o{v}_{\*}(z_2,s)),\I^{\p}(\o{v}_{\*}(z_2,s))\times\\
	&&(p+2\delta z_1),
	-\I^{\p}(\o{v}_{\*}(z_2,s))X_2+\I^{\p\p}(\o{v}_{\*}(z_2,s))p\otimes p)\\
	&& \leq \frac{1}{\sqrt{\l_0}}\nu_{1,R_{\alpha}}(|\I^{\p}(\o{v}_{\*}(z_2,s))(p+2\delta z_1)-
	\I^{\p}(\o{v}_{\*}(z_2,s))p|)\\
	&& \leq \frac{1}{\sqrt{\l_0}}\nu_{1,R_{\alpha}} (2\sqrt{\delta\L_0}),
\end{eqnarray*} 
since that 
\[|\I^{\p}(\o{v}_{\*}(z_2,s))(p+2\delta z_1)|,\;\;\;|\I^{\p}(\o{v}_{\*}(z_2,s))p|\leq R_{\alpha}\]
and
\[|-\I^{\p}(\o{v}_{\*}(z_2,s))X_2+\I^{\p\p}(\o{v}_{\*}(z_2,s))p\otimes p|\leq R_{\alpha}\]

The inequality (\ref{CP5}) is satisfied for
\begin{eqnarray*}
X&=&\I^{\p}(\o{v}_{\*}(z_2,s))(X_1+2\delta I_N)+2\delta\I^{\p\p}(\o{v}_{\*}(z_2,s))(p\otimes z_1+z_1\otimes p)\\
&&+4\delta^2\I^{\p\p}(\o{v}_{\*}(z_1,s))z_1\otimes z_1
\end{eqnarray*} 
\[Y=\I^{\p}(\o{v}_{\*}(z_2,s))X_2.\] 

\[Z=\I^{\p\p}(\o{v}_{\*}(z_2,s))p\otimes p\]
In fact, (\ref{CP12}) implies
\begin{eqnarray}\label{CP14}
		&& -\varepsilon_1 I_{2N}\leq\left(\begin{array}{cc}
					X & O\\
					O & Y
			\end{array}	\right)=\I^{\p}(\o{v}_{\*}(z_2,s))   
                                          \left(\begin{array}{cc}
					X_1 & 0\\
				        0 &  X_2
			\end{array}	\right)\nonumber \\ &&+
					\left(\begin{array}{cc}
					2\delta \I^{\p}(\o{v}_{\*}(z_2,s))I_N & 0\\
				        0 &  0
			\end{array}	\right)\nonumber\\ 
					&&+2\delta\alpha \I^{\p\p}(\o{v}_{\*}(z_2,s)) 
					\left(\begin{array}{cc}
					(z_1-z_2)\otimes z_1+z_1\otimes(z_1-z_2) & 0\\
				        0 &  0
			\end{array}	\right)\nonumber\\ 
					&&+4\delta^2 \I^{\p\p}(\o{v}_{\*}(z_2,s)) 
					\left(\begin{array}{cc}
					z_1\otimes z_1 & 0\\
				        0 &  0
			\end{array}	\right)\nonumber\\ 
					&&\leq 3\alpha\sqrt{\L_0}
				\left(\begin{array}{cc}
					I_N  & -I_N\\
					-I_N & I_N
			\end{array}	\right)+ \nonumber\\ 
				&&+2\big(\delta\sqrt{\L_0} + (\alpha \sqrt{\delta}+\delta) \|z^{\p}\|_{\infty}\big)
              				\left(\begin{array}{cc}
					 I_N & 0\\
				        0 &  I_N
			\end{array}	\right).
\end{eqnarray}
So the relation (\ref{CP5}) holds if we choose $\e_2=3\alpha\sqrt{\L_0}$, 
$\e_3=2\big(\delta\sqrt{\L_0} + (\alpha \sqrt{\delta}+\delta) \|z^{\p}\|_{\infty}\big)$, $\e_1=\e_2+\e_3$. 
If we choose $R$ as in (\ref{CP5}), indipendent of $\delta$, and $R\geq R_{\alpha}$, then holds,
\begin{eqnarray*}
	&&\frac{1}{\I^{\p}(\o{v}_{\*}(z_2,s))}\Big[F(z_2,s,\I(\o{v}_{\*}(z_2,s)),\I^{\p}(\o{v}_{\*}(z_2,s))(p+2\delta z_1),
	-\I^{\p}(\o{v}_{\*}(z_2,s))X_2+\\
	&&\I^{\p\p}(\o{v}_{\*}(z_2,s))p\otimes p)\\
	&&-F(z_1,s,\I(\o{v}_{\*}(z_2,s)),\I^{\p}(\o{v}_{\*}(z_2,s))(p+2\delta z_1),
	\I^{\p}(\o{v}_{\*}(z_2,s))(X_1+2\delta I_N)+\\
	&&\I^{\p\p}(\o{v}_{\*}(z_2,s))(p+2\delta z_1)\otimes(p+2\delta z_1))\Big]\leq\\
	&&\frac{1}{\sqrt{\l_0}}\nu_{2}\big(|z_1-z_2|(1+\sqrt{\L_0}\alpha|z_1-z_2|+2\sqrt{\L_0\delta})+
	\e_2 |z_1-z_2|^2 \big)\\
&&+\frac{1}{\sqrt{\l_0}}\nu_{2,R}(2 \e_3).
\end{eqnarray*}
If we take, in (\ref{CP7}) $\l=\I^{\p}(\u{v}^{\*}(z_1,s))$, $\w{\l}=\I^{\p}(\o{v}_{\*}(z_2,s))$,
$\kappa=\I^{\p\p}(\u{v}^{\*}(z_1,s))$, $\w{\kappa}=\I^{\p\p}(\o{v}_{\*}(z_2,s))$, $q=p+2\delta z_1$,
$X=X_1+2\delta I_N$ and 
\begin{eqnarray*}
	&&\frac{1}{\I^{\p}(\o{v}_{\*}(z_2,s))}F(z_1,s,\I(\o{v}_{\*}(z_2,s)),\I^{\p}(\o{v}_{\*}(z_2,s))(p+2\delta z_1),
\I^{\p}(\o{v}_{\*}(z_2,s))\times\\\\
&&(X_1+2\delta I_N)+
\I^{\p\p}(\o{v}_{\*}(z_2,s))(p+2\delta z_1)\otimes(p+2\delta z_1))-
\frac{1}{\I^{\p}(\u{v}^{\*}(z_1,s))}\times \\\\ && F(z_1,s,\I(\u{v}^{\*}(z_1,s)),\I^{\p}(\u{v}^{\*}(z_1,s))(p+2\delta z_1),
\I^{\p}(\u{v}^{\*}(z_1,s))\times (X_1+2\delta I_N)+\\\\
&&\I^{\p\p}(\u{v}^{\*}(z_1,s))(p+2\delta z_1)\otimes(p+2\delta z_1))
\leq \Gamma_0(\u{v}^{\*}(z_1,s)-\o{v}_{\*}(z_2,s))
\end{eqnarray*}
In fact $\I^{\p}(v)= \sqrt{z(\I(v))} \in [\sqrt{\l_0},\sqrt{\L_0}]$ and $\I^{\p\p}(v)=\textstyle{\frac{z^{\p}(\I(v))}{2}}$
Replacing the obtained estimates in the inequality (\ref{CP15}), yields
\begin{eqnarray*}
	\varphi^{\p}(s)&\leq&  \nu_{1,R_{\alpha}} (2\sqrt{\delta\L_0})+ \\
	&&+\frac{1}{\sqrt{\l_0}}\nu_{2}
	\big(|z_1-z_2|(1+\sqrt{\L_0}\alpha|z_1-z_2|+2\sqrt{\L_0\delta})+
	\e_2 |z_1-z_2|^2 \big)+\\
	&&+\nu_{2,R}(2 \e_3)+\Gamma_0(\u{v}^{\*}(z_1,s)-\o{v}_{\*}(z_2,s)).
\end{eqnarray*} 
Letting $\delta\rightarrow 0$ (without moving $\alpha $), and considering the assertion $(iv)$ of Lemma \ref{LEMMATECNICO}, yields,
\begin{eqnarray*}
	\varphi^{\p}(\lim_{\delta\rightarrow 0}s)&\leq& 
	\Gamma_0(\lim_{\delta\rightarrow 0}\u{v}^{\*}(z_1,s)-\o{v}_{\*}(z_1,s))\\&&+
	\frac{1}{\sqrt{\l_0}}\nu_{2}
	\big(\lim_{\delta\rightarrow 0} (|z_1-z_2|+4\sqrt{\L_0}\alpha|z_1-z_2|^2))
\end{eqnarray*} 
then, again with condition $(iv)$ and $(vi)$ of Lemma \ref{LEMMATECNICO},
 letting $\alpha\rightarrow \infty$ in the above inequality, yields,
\begin{eqnarray*}
	\varphi^{\p}(t_0)\leq \Gamma_0 (\vartheta^{\*}(t_0))
\end{eqnarray*}
By Lemma \ref{OSGOOD}, and the monotonicity of the application $\Psi$, the comparison assertion
(\ref{CP16}) holds.
\cvd
\end{Dimo}


In order to show why our conditions on the Hamiltonian $F$ are more general with respect to
the hypothesis up to now proposed in the viscosity theory to state the comparison principle,
we build some examples which do not satisfy usual assumptions for the standard comparison result,
with particular regard to the recent paper of P.L Lions and P.E. Souganidis \cite{LS}.\\
In the example \ref{example1}, we consider an Hamiltonian which has not the Lipschitz regularity with respect to $u$. In that case
the appeal to the Osgood property becomes a necessary requirement to obtain the comparison. For seeing that our
structural approach is quite different than the technical assumptions proposed in \cite{LS} and 
following the same kind of setting also used in \cite{LS} for solving a stochastic pde, in the example \ref{example2} we give an Hamiltonian
$F_{\Phi}$, which comes from a function $F=F(\nabla u,\nabla^2 u)$, through a regular transformation $\Phi$. The first function $F$, which we propose in such example, 
is not Lipschitz continuous, as instead require the authors of \cite{LS} in order to apply their techniques. Our method not only, seems to be efficient when there 
is a lack of regularity in the Hamiltonian, but also when we lose the hypothesis (1.12), pag. 621 in \cite{LS}, see (\ref{LS2}) below. Actually, 
in example (\ref{example2}) we consider another function $F$ which is Lipschitz continuous, but does not satisfy that condition. In both cases 
our Theorem \ref{CP} can be successfully applied. Finally, before to begin the study of the financial model (\ref{DM1}),
in the observation \ref{MBSmodel}, we remark that it seems difficult to directly combine the results 
of \cite{LS} for deducing the comparison principle for the equation (\ref{DM1}). 

\begin{example}\label{example1}
Let $F:(-1,\infty)\times\R^N\times\S^N$ be the Hamiltonian defined by
\begin{eqnarray}\label{F1}
F(u,p,X)&=&-tr(X)+\frac{1}{u+1}|p|^2+\varphi(u),\nonumber\\
\varphi(u)&=&\left\{\begin{array}{cc}
	   	(u^2+u)\log(u) & \mbox{if $u>0$,}\\
		0 & \mbox{if $u\leq 0$.}	
\end{array}\right.
\end{eqnarray}
This function does not get back in the class of classical problem, for which the comparison 
is already stated (see \cite{1}, \cite{6}); actually we have a quadratic gradient bounds and also a lack of regularity in $u$.
Of course $F$ satisfies conditions (\ref{CP2}), (\ref{CP3}) of Theorem \ref{CP} in its domain of definition, while
in the interval $I\equiv[-\frac{1}{2},\frac{1}{e}]$ verifies (\ref{CP4})-(\ref{CP6}). By the regularity of $F$ with respect to $p$ and $X$, to verifying
condition (\ref{CP7}), it suffices to prove the existence of a function $z=z(u)>0$ such that,
\begin{eqnarray}\label{F2}
F_z(u,p,X)&=&\frac{1}{\sqrt{z(u)}}F(u,\sqrt{z(u)}p,\sqrt{z(u)}X+\frac{z^{\p}(u)}{2}p\otimes p)\nonumber\\
&=& -tr(X)+\big(\frac{\sqrt{z(u)}}{u+1}-\frac{z^{\p}(u)}{2\sqrt{z(u)}}\big)|p|^2+\frac{\varphi(u)}{\sqrt{z(u)}},
\end{eqnarray}
satisfies
\begin{eqnarray}\label{F3}
&&F_z(u,p,X)-F_z(v,p,X)\geq -\Gamma(u-v),\nonumber\\
&&\;\;\;\;\;\;\;\forall\;-\frac{1}{2}\leq v< u\leq \frac{1}{e},\;(p,X)\in\R^N\times\S^N,
\end{eqnarray}
with $\Gamma$ an Osgood-type function. By the previous consideration about the dependence on $u$,
we cannot turn to the $u$-partial derivative of $F_z$, to obtain estimate (\ref{F3}). 
Consider the function
\begin{eqnarray}\label{F4}
\Gamma(h)=\left\{\begin{array}{cc} 
		h\log(\frac{1}{h}) & \mbox{if $0<h<\frac{1}{e}$},\\
		& \\
		\frac{1}{e} & \;\;\;\;\;\;\;\mbox{if $\frac{1}{e}\leq h\leq \frac{1}{e}+\frac{1}{2}$},
\end{array}\right.
\end{eqnarray}
for $h>0$; then $\Gamma$ is an Osgood-type function in the sense of definition \ref{DOSGOOD}, while it is a simple
exercise to show that $\Gamma(h)=\sup_{x\in I}[x\log(x)-(x+h)\log(x+h)]$. Therefore, choosing $z(u)=(u+1)^2$, the inequality (\ref{F3}), becomes
\begin{eqnarray}\label{F5}
u\log(u)-v\log(v)\geq
-\Gamma(u-v),
\end{eqnarray}
which by the definition of $\Gamma$ it holds. Hence, applying Theorem \ref{CP}, we have the comparison between sub/super solutions
which take values in $I$, for (\ref{F1}).
\end{example}

\begin{oss}\label{LS1}\rm
In \cite{LS}, it was proposed a comparison principle for a class of problems which do not satisfy the usual assumptions.
To treat a stochastic differential problem the authors have to prove a comparison between bounded solutions, for a deterministic problem.
They consider the following equation,
\begin{eqnarray}\label{GPT}
	\partial_t v &=& F_{\Phi}(v,\nabla v,\nabla^2 v), \nonumber\\
	F_{\Phi}(u,p,X) &=& \frac{1}{\Phi^{\p}(u)}F(\Phi^{\p}(u)p,\Phi^{\p}(u)X+\Phi^{\p\p}(u)p\otimes p),\nonumber\\
	&&\;\;\;\forall\;\;(u,p,X)\in\R\times\R^N\times\S^N,
\end{eqnarray}
in $\R^N\times (0,T)$, where $\Phi$ is a smooth function, $F$ is independent of the unknown and
$F(\cdot,X)\leq F(\cdot, Y)$, if $X\leq Y$. Moreover the function $F$ is a globally Lipschitz continuous function of their variables, and
satisfies a structural condition (see (1.12), pag. 621 in \cite{LS}),
\begin{eqnarray}\label{LS2}
\left\{\begin{array}{c}
X\cdot \nabla_X F+p\cdot\nabla_p F-F\leq C\;\;\; \mbox{ for a.e. $(X,p)$}\\
or \\
X\cdot \nabla_X F+p\cdot\nabla_p F-F\geq C\;\;\; \mbox{ for a.e. $(X,p)$}\\
\end{array}\right.
\end{eqnarray}
for some constant $C>0$. Under these assumptions, making a global change of the unknown $v=\phi(w)$, by a transformation $\phi$, they prove that
the partial derivative of the new Hamiltonian with respect to the unknown $w$ is bounded from above.
This fact, implies that the comparison result follows from the classical theory.\\
If we consider the problem (\ref{GPT}), coming from a Lipschitz continuous function $F$ which satisfies (\ref{LS2}), we
can define $z(v)=[\phi^{\p}(\phi^{-1}(v))]^2$ to verify condition (\ref{CP7}) of Theorem \ref{CP}; while the other
properties required by Theorem \ref{CP} are easily derived by the degenerate ellipticity, and the regularity assumptions on $F$, and $\Phi$.
Considering the previous setting we build two Hamiltonians $F_{\Phi}$, which 
satisfy our conditions but such that it is not possible to use the methods of \cite{LS}.
\begin{example}\label{example2}
In many situations the Hamiltonian $F$, is not Lipschitz continuous, so the conditions (\ref{LS2}), cannot be verified. 
Let $\gamma$ be a number in $(0,1)$, and consider the 1-dimensional Hamiltonian
\begin{eqnarray}\label{LS5}
\t{F}(v,q,Y)=Y-\frac{1}{v}q^2-v^{1-\gamma}|q|^{\gamma},\;\;\;\forall\;(v,q,Y)\in(0,\infty)\times\R\times\R.
\end{eqnarray}
Of course $\partial_v\o{F}$ is not bounded from above. Moreover if $\Phi(v)=\log(v)$, and $F(p,X)=X-|p|^{\gamma}$, then 
\begin{eqnarray}\label{LS6}
\t{F}=F_{\Phi}.
\end{eqnarray}
The assumptions (\ref{CP3}), (\ref{CP4}), (\ref{CP6}) follow by the linearity of $F$ with respect to $X$ and the
holder continuity of $F$ with respect to $p$. Now it is an exercise to prove that choosing $z(v)=e^{-2v}$, we obtain
also condition (\ref{CP7}). Hence although the lack of regularity for $F$, we could apply Theorem \ref{CP}.\\
As a second example consider,
\begin{eqnarray}\label{LS7}
F(X,p)&=& X+g(p),\;\;\;(p,X)\in\R^2\nonumber\\
g(p)&=& \left\{\begin{array}{cc}
	\;\;\,\log(1+p) & \mbox{if $p>0$}\\
	-\log(1-p) & \mbox{if $p\leq 0$}
\end{array}\right. \nonumber\\
\Phi(u) &=& \arctan(u).
\end{eqnarray}
Then $F$ is a Lipschitz continuous function in $\R^2$, and $\Phi$, has the regularity required in \cite{LS}.
Nevertheless $p\cdot g^{\p}(p)-g(p)$, is neither bounded from above nor from below, so (\ref{LS2}) does not hold. 
To applying our procedure, as in example \ref{example1} it suffices to prove that the Hamiltonian
\begin{eqnarray}\label{LS8}
(F_{\Phi})_z(u,p,X) &=& X+l(u)p^2+\frac{g(w(u)p)}{w(u)},\nonumber\\
l &=& (\sqrt{z})^{\p}+\frac{\Phi^{\p\p}}{\Phi^{\p}}\sqrt{z},\nonumber\\
w &=& \Phi^{\p}\sqrt{z},
\end{eqnarray}
satisfies $\Delta\equiv(F_{\Phi})_z(u,p,X)-(F_{\Phi})_z(v,p,X)\leq \gamma\cdot (u-v)$, for an appropriate function $z$, $u\geq v$ in a bounded interval, and $\gamma>0$.\\
Taking $z(u)=(u^2+1)^2 \cdot (\beta-(1/2)\arctan^2(u))^2$, with $8 \beta>\pi^2$ we have
$l(u)=-\arctan(u)$, $w(u)=\beta-(1/2)\arctan^2(u)$. Then by the definition of $g$ and (\ref{LS8}), 
\begin{eqnarray}\label{LS9}
\Delta\leq \frac{1}{\t{w}^2}\big[\big(l(u)-l(v)\big)|\t{w}p|^2+\alpha\big(1+\log(1+|\t{w}p|)\big)\big],
\end{eqnarray}
where $\t{w}>0$, is a point between $w(u)$ and $w(v)$, $\alpha>0$ depends on the bounds of $u$. Therefore,
by $l^{\p}<0$, it follows the assertion.
\end{example}
\end{oss}
The previous observation and the last two examples show that our com\-pa\-ri\-son principle really extend the result of \cite{LS}.

\begin{oss}\label{MBSmodel}\rm
As we have already said in the introduction of the present paper, we can neither  apply the 
the existing comparison results nor the same technique used in \cite{LS}, to obtain the result for the financial model.
In the classical theory was not considered a quadratic dependence in the variable $p$. In \cite{LS}
the extension to problems with this kind of growth was examined, but even if we eliminate the dependence
by the spatial variable $x$ in (\ref{DM1}), changing in this way the structure of the model, we can reduce the equation 
to a pde, whose Hamiltonian has the form (\ref{GPT}), where $F(p,X)=-(1/2)tr(\s\s^{\top}X)+C|\s^{\top}p+f|^2 $, for some constants $C$ and $f=(f_1,\ldots,f_d)$. 
Nevertheless this Hamiltonian is not globally Lipschitz continuous.
\end{oss}


\section{Application to the Financial Model}\label{MBS}
In this section we give an application of Theorem \ref{CP}, for treating a real
problem which comes from the mathematical finance, and which is described by the equation (\ref{DM1}).
We shall make the following assumptions:
\begin{description}
\item[(${\bf P_1}$)]   $\mu$ is bounded, and $\mu(\cdot,t)\in Lip(\R^N)$, for every $t\in [0,T)$,
		       with a Lipschitz constant indipendent of the time.
\item[(${\bf P_2}$)] $\rho,\;\tau>0$ are known real parameters, 
		       $h\in {\cal C}^{2,1}(\R^{N}\times[0,T))$, is bounded from above and nonnegative. Moreover 
		       $\nabla h$ is bounded, and, 
                       $\partial_t h(\cdot,t)$, $tr(\sigma\sigma^{\top}(t) \nabla^2 h(\cdot,t)),\;\nabla h(\cdot,t)\in Lip(\R^N)$ 
		       with Lipschitz constants indipendent of the time $t\in[0,T)$. $\xi\in \C^1([0,T])$, is positive.
\item[(${\bf P_{3}}$)] $U_0\in Lip(\R^N)$, with $U_0\geq 0$, and bounded from above.  			   		                                       
\end{description}


\begin{oss}
The assumption about the sign of $h$ and $\xi$, reflect the financial modelization:
$h$ describes the cash-flow from the pool of mortgages, while $\xi$ represents the value of some
Bank-account, which moves with the interest rate $r$. For more generality, the initial datum $U_0$ is assumed
nonnegative and regular because, in a real framework, it is constant
and equal to zero. Actually at the maturity date $T$, the mortgage-holders have
paid off their debt.\\
The constants $\rho,\;\tau>0$ have a financial interpretation, in fact
they represent, respectively, the risk-aversion coefficient and the coupon rate paid by the mortgage-holder.
\end{oss}

Using Theorem \ref{CP}, it is now easily to show the comparison for equation (\ref{DM1}).

\begin{teo}\label{DM1COMPARISON}
Assume conditions (${\bf P_1}$), (${\bf P_2}$), (${\bf P_3}$). Set the functions
\begin{eqnarray}\label{k}
\u{k}(t)&:=&\!\!\! e^{-\int_0^t r(x)dx}\Big( \inf_{\R^N} U_0+\int_0^t e^{\int_0^s r(x)dx}\inf_{\R^N}[(\tau-r(s))h(x,s)]ds\Big)\nonumber\\
\o{k}(t)&:=&\!\! K_0 t+c_0,\;\;\;\;\;\;\;c_0:=\max\Big(\sup_{\R^N}U_0,\sup_{[0,T)}\u{k}\Big),\nonumber\\
 && \!\!K_0:=\!\!\sup_{\R^N\times [0,T)}\Big(\frac{(\tau-r(t))h(x,t)-c_0 r(t)}{1+t\cdot r(t)}\Big)^{+},\;\;\!\!t\in[0,T). 
\end{eqnarray}
If exists a modulus $\nu$ such that 
\begin{eqnarray}\label{d_t h}
|h(x,t)-h(x,s)|\leq \nu(|t-s|),
\end{eqnarray}
and
\begin{eqnarray}\label{XI}
\xi(t)+h(x,t)+\u{k}(t)>0,
\end{eqnarray}
for every $t,\;s\in [0,T)$, $x\in\R^N$, then $\u{k}\leq \o{k}$ are, respectively, $\C^1([0,T))$ sub/supersolution of (\ref{DM1}),
and for every sub/supersolution $\u{U},\;\o{U}$ of (\ref{DM1}), such that $\u{k}\leq \u{U},\;\o{U}\leq \o{k}$, then
$\u{U}\leq \o{U}$ holds in $\R^N\times [0,T)$. 
\end{teo}

\begin{Dim}. Since the inequalities (\ref{d_t h}), (\ref{XI}), the function $\u{k}$, defined by (\ref{k}), is a regular subsolution of (\ref{DM1}):
actually (\ref{d_t h}) implies continuity of $\inf_{x\in\R^N}[(\tau-r(t))h(x,t)]$ as a function of the time. It
is also easy to see, by the definition of $K_0$, that $\o{k}$ is a supersolution of the same problem.
By assumptions (${\bf P_1}$)-(${\bf P_3}$), we can consider the change $u=U+h+\xi$. Then
if $U$ is a vi\-sco\-si\-ty sub/supersolution of (\ref{DM1}) with $\u{k}\leq U\leq \o{k}$, then by (\ref{XI}),
$u$ is a positive viscosity sub/supersolution of
\begin{eqnarray}\label{DM2}
\!\!\partial_t u\!-\!\frac{1}{2}tr(\s\s^{\top}\nabla^2 u)\!-\!\la \mu,\nabla u \ra\!+\!\rho\frac{|\s^{\top}\nabla u\!-\!\s^{\top}\nabla h|^2}{u}
\!+\!r(t)u\!+\!g(x,t)=0,
\end{eqnarray}
in $\R^N\times(0,T)$, where $g:=-\partial_t h+\frac{1}{2}tr(\s\s^{\top}\nabla^2 h)+\la \mu,\nabla h\ra
-\tau h -(\xi^{\p}+r(t)\xi)$. Moreover
by the conditions on $h$, the new initial datum $u_0=U_0+h(x,0)+\xi(0)$, is globally
Lipschitz continuous. For proving the comparison through Theorem \ref{CP} it suffices to show that conditions
(\ref{CP2})-(\ref{CP7}) hold for equation (\ref{DM2}), when we consider the interval $[a,b]=[\inf(\u{k}+h+\xi), K_0 T+c_0+\sup (h+\xi)]$. 
Conditions (\ref{CP2})-(\ref{CP4}) follow
immediately, while conditions (\ref{CP5})-(\ref{CP6}) follow by the linearity of the second order part in the equation (\ref{DM2})
and by the Lipschitz regularity for $g$ and assumptions (${\bf P_1}$), (${\bf P_2}$). 
Set $m_0:=\inf_{\R^N\times [0,T]}(\u{k}+h+\xi)$, $M_0:=K_0 T+c_0+\sup(h+\xi)$ and consider the function $z(u)=(\l_1 u-\l_2)^2$, $u\in [m_0, M_0 ]$, for constants $\l_1,\;\l_2$ such that $\l_1 m_0>\l_2>0$. Now we verify
(\ref{CP7}) with this choise. Let $\l,\;\w\l,\;\kappa,\;\w\kappa$ be as in condition (\ref{CP7}) of Theorem \ref{CP}, $(x,t,X)\in\R^N\times [0,T)\times\S^N$, $R>0$ and $|q|\leq R$, then for 
$M_0\geq u>v\geq m_0$, we have,
\begin{eqnarray*}
\lefteqn{\frac{1}{\l}\big[-\frac{1}{2}tr(\s\s^{\top}(\l X\!\!+\!\! \kappa q\otimes q))-\la \mu,\l q \ra+\rho\frac{|\s^{\top}\l q-\s^{\top}
\nabla h|^2}{u}+r(t)u\!+\!g(x,t)\big]-}\nonumber\\
&&\frac{1}{\w\l}\big[\!\!-\frac{1}{2}tr(\s\s^{\top}(\w\l X\!\!+\!\!\w\kappa q\otimes q))\!-\!\la \mu,\w\l q \ra\!+\!\rho\frac{|\s^{\top}\w\l 
q-\!\!\s^{\top}\nabla h|^2}{v}+r(t)u\!\!+\!\!g(x,t)\big]\nonumber\\
&\geq&-\frac{z^{\p}(u)}{4\l}|\s^{\top}q|^2+\rho\frac{ |\l \s^{\top}q-\s^{\top}\nabla h|^2}{\l u}+
\frac{z^{\p}(v)}{4\w\l}|\s^{\top}q|^2-\rho\frac{|\w\l\s^{\top}q-\s^{\top}\nabla h|^2}{\w\l v}\nonumber\\
&\geq &\rho\big[\frac{\sqrt{z(u)}}{u}-\frac{\sqrt{z(v)}}{v}\big]|\s^{\top}q|^2+
2\frac{\rho}{v}\la\s^{\top}q,\s^{\top}\nabla h\ra-2\frac{\rho}{u}\la\s^{\top}q,\s^{\top}\nabla h\ra
+\nonumber\\ &&\rho|\s^{\top}\nabla h|^2 
\big(\frac{1}{\sqrt{z(u)}u}-\frac{1}{\sqrt{z(v)}v}\big)
-C_1\big((|z(u)-\l^2|+|z(v)-\w\l^2|)(1+|q|)\big)\nonumber\\
&\geq&\rho\big(\frac{\l_1 u-\l_2}{u}- 
\frac{\l_1 v-\l_2}{v}\big)|\s^{\top}q|^2-\rho C_2(|\s^{\top}q|+1)(u-v)\nonumber\\
&&-C_1\big((|z(u)-\l^2|+|z(v)-\w\l^2|)(1+|q|)\big)\nonumber\\
&\geq&\rho\frac{\l_2}{M_0^2} (u-v)|\s^{\top}q|^2-\rho C_2(|\s^{\top}q|+1)(u-v)
-C_1\big((|z(u)-\l^2|\nonumber\\
&&+|z(v)-\w\l^2|)(1+|q|)\big)\geq
\end{eqnarray*}
\begin{eqnarray}\label{VERIFICATION}
&\geq& -\rho(u-v)(\frac{(C_2 M_0 )^2}{4\l_2}+C_2)-C_1\big((|z(u)-\l^2|\nonumber\\
&&+|z(v)-\w\l^2|)(1+|q|)\big),
\end{eqnarray}
where, \[C_1=C_1(R,\l_1,\l_2,m_0, M_0, \|\s^{\top}\|_{\infty},\|\nabla h\|_{\infty})\] \[C_2=C_2(\l_1,\l_2,m_0, M_0,
\|\s^{\top}\|_{\infty},\|\nabla h\|_{\infty})\] are positive constants. 
By (\ref{VERIFICATION})
we see that $\Gamma(x)=\rho(\frac{(C_2 M_0 )^2}{4\l_2}+C_2)x$ is a function as in Definition
\ref{DOSGOOD}, and conditions of Theorem \ref{CP} are satisfied for problem (\ref{DM2}). The proof of
comparison is now a direct consequence of Theorem \ref{CP} applied to (\ref{DM2}).\cvd
\end{Dim}

\begin{oss}
Also the functions $\u{k},\;\o{k}$ defined in (\ref{k}), have a financial interpretation. Actually
in the real differential problem for $MBS$, where $U_0\equiv 0$, and after the usual time-change of variable $T-t=s$, they represent, respectively,
the evolution of a Zero-Coupon-Bond with maturity $T$, which follows the movements of $r$, and
a linear above estimate in the time for the evolution of the remaining principal in the pool.  
\end{oss}

\begin{cor}\label{Fine}
Under the same assumptions of Theorem \ref{DM1COMPARISON}, the problem
(\ref{DM1}) has a unique continuous viscosity solution $U$, such that
$\u{k}\leq U\leq \o{k}$ over $\R^N\times [0,T)$.
\end{cor}

\begin{Dim}. By Theorem \ref{DM1COMPARISON}, we have proved the comparison between viscosity sub\-su\-per\-so\-lu\-tion
of problem (\ref{DM1COMPARISON}), then the assertion follows adapting a version of Per\-ron's Me\-thod for
the elliptic case by H. Ishii in \cite{21} to the parabolic case.\cvd
\end{Dim}

Even though, in (\ref{DM2}), we have a term which depends on the unknown $u$ which multiplies a one-order term
in a descreasing form with respect to the variable $u$, as is showed in the following Proposition \ref{approx1}
is possible to prove a regularity result of the viscosity solution. 
As we prove in the following Theorem, we can consider a general class of quasilinear
equations for which is possible to obtain a Lipschitz regularity for the viscosity solution.

\begin{teo}\label{approx1}
Let $v$ be a continuous viscosity solution of the following problem
\begin{eqnarray}\label{approx2*}
&&\partial_t v-\frac{1}{2}tr(\sigma\sigma^{\top}\nabla^2 v)
-\langle \mu,\nabla v\rangle+\l_1(v)                
|\sigma^{\top}\nabla v|^2+\l_2(v)\la\sigma^{\top} \nabla v,w\ra \nonumber\\\nonumber\\
&&+f(x,t,v)=0,\;\;\;\;\;\;\;(x,t)\in\R^N\times(0,T).
\end{eqnarray}
With an initial datum $v_0\in Lip(\R^N)$. Suppose that $v$ takes values in the interval $[c,d]$,
$\l_1\lambda_2\in \C^1([c,d])$, 
$\l_1^{\p}>0$ on $[c,d]$, $\s$ $\mu$ are as in problem (\ref{DM1}), with $\mu,\;w$ satisfying assumption \em{(}${\bf P_1}$\em{)}, 
$f\in C(\R^N\times[0,T)\times [c,d])$ and
$f(\cdot,t,\cdot)\in Lip(\R^N\times [c,d])$ with a Lipschitz constant wich is independent of the time. 
If $M>\frac{\textstyle{Lip(v_0)}}{\textstyle{2}}$, and 
\begin{eqnarray}\label{Constant}
\lefteqn{C=2Lip(\mu)+\frac{\|\l^{\p}_2\|^2_{\infty}\|w\|^2_{\infty}}{4\min_{[c,d]}\l_1^{\p}}+}\nonumber\\
&&+2\|\l_2 \|_{\infty} \|{\sigma}^{\top} \|_{\infty} Lip(w)+Lip(f)(\frac{1}{2 M}+1)
\end{eqnarray}
then 
\begin{eqnarray}\label{approx20}
|v(x,t)-v(y,t)|\leq 2 M e^{C t}|x-y|
\end{eqnarray}
holds for every $x,\;y\in\R^N $, $t\in [0,T)$. In particular for every $t\in [0,T)$, $v(\cdot,t )\in Lip(\R^N)$.
\end{teo}  

As consequence of Theorem \ref{approx1}, we have a regularity result for our differential model.

\begin{prop}
Assume the same assumptions of Theorem \ref{DM1COMPARISON}.
If $U$ is the viscosity solution of (\ref{DM1}), then $U(\cdot,t)\in Lip(\R^N)$, with a Lipschitz constant which
is independent of the time.
\end{prop}

First use Theorem \ref{approx1} to prove this Proposition, then we shall give the proof of the Theorem.\\

\begin{Dim}.
Consider $m_0,\;M_0$ and the change $U+h+\xi=u$, already used in the proof of Theorem \ref{DM1COMPARISON}. Then
$u$ solves, in a viscosity sense the equation (\ref{DM2}), with an initial datum $u_0\in Lip(\R^N)$; moreover, by assumption
(${\bf P_3}$),  $u\in [m_0,M_0]$. By the regularity of $h$, it suffices to prove the assertion for the function $u$. For a regular trasformation $\I$, defined
in a open neighbourhood of some interval $[c,d]$, with $\I^{\p}>0$, $\I([c,d])=[m_0,M_0]$, the function $v:=\I^{-1}(u)$, is a solution
in a viscosity sense of a Cauchy problem which has the same structure of problem (\ref{approx2*}) with,
\begin{eqnarray}\label{approx21}
\l_1(v) &=&\frac{d}{d v}\log\Big(\frac{\I^{\rho}(v)}{\sqrt{\I^{\p}(v)}}\Big),\;\;
\l_2(v) =-\frac{2\rho}{\I(v)},\\
w(x,t)&=& \s^{\top}(t)\nabla h(x,t),\\
f(x,t,v)&=& \frac{|\s^{\top}(t)\nabla h(x,t)|^2}{\I(v)\I^{\p}(v)}+\frac{g(x,t)}{\I^{\p}(v)}
+r(t)\frac{\I(v)}{\I^{\p}(v)},
\end{eqnarray}
for every $(x,t,v) \in \R^N \times [0,T)\times [c,d]$. Morevover the initial datum is $v_0=\I^{-1}(u_0)$. 
For defining the transformation $\I$, we take the same kind of function $z$, which we have used for proving Theorem \ref{DM1COMPARISON},
\begin{eqnarray}\label{approx22}
\I(v)=m_0\frac{e^{\frac{ 2 v}{m_0}}+1}{2},\;\;\;\;\forall\;\;v\in\Big[0,\frac{m_0}{2}\log\Big(\frac{2 M_0}{m_0}-1\Big)\Big].
\end{eqnarray}
This  choise, yields,
\begin{eqnarray}\label{approx23}
\l_1^{\p}(v)=\frac{4\rho}{m_0^2}\frac{e^{\frac{2 v}{m_0}}}{(e^{\frac{2 v}{m_0}}+1)^2}>0.
\end{eqnarray}
By the assumption (${\bf P_1}$), (${\bf P_2}$) it follows that $w$, $f$ satisfy the conditions of Theorem \ref{approx1}, so applying this Theorem
we infer that exists a constant $C>0$, such that 
\begin{eqnarray}\label{approx24}
|v(x,t)-v(y,t)|\leq M e^{C t}|x-y|,\;\;\;\;\forall\;x,\;y\in\R^N,\;t\in [0,T),
\end{eqnarray}
with $M>Lip(v_0)$, and $C=C(M)$; hence the definition (\ref{approx22}), yields
\begin{eqnarray}\label{approx25}
|u(x,t)-u(y,t)|\leq  M \Big(\frac{2 M_0}{m_0}-1\Big) e^{C t}|x-y|,\;\;\forall\;x,\;y\in\R^N,\;t\in[0,T).
\end{eqnarray}
\cvd
\end{Dim}

\begin{Dimo} {\bf{Theorem \ref{approx1}.}}
Consider the function $Q$ defined by $Q:=v e^{-Ct}$, where $C\geq 0$ is a positive constant given by (\ref{Constant})
the function $Q$ is a continuous viscosity solution of the equation
\begin{eqnarray}\label{approx8}
&&\partial_t Q-\frac{1}{2}tr(\s \s^{\top} \nabla^2 Q)-\la \mu,\nabla Q\ra+e^{Ct}\l_1(Q e^{Ct})|\s^{\top}\nabla Q|^2
\nonumber\\
&&+\l_2(Q e^{Ct})\la \s^{\top} \nabla Q, w \ra+e^{-Ct} f(x,y,t,Q e^{Ct})+CQ=0.
\end{eqnarray}
With the same initial datum $v_0$. For $\gamma,\;\delta,\;\var>0$, $2M>Lip(v_0)$ we set
\begin{eqnarray}\label{approx2}
 H(x,y,t):= Q(x,t)-Q(y,t)-\var|x|^2-K(x,y,t),\nonumber\\
 K(x,y,t):= M \big( \frac{|x-y|^2}{\delta}+\delta \big)+\frac{\gamma}{T-t},\;\;\;\;(x,y,t)\in\R^N\times[0,T).
\end{eqnarray}
We will show that for every $\delta,\;\gamma>0$  there is $\var_0=\var_0(\delta,\gamma)>0$ such that
for $0<\var<\var_0$ we have,
\begin{eqnarray}\label{approx19}
H(x,y,t)\leq 0,\;\;\forall\;(x,y,t)\in\R^N\times[0,T). 
\end{eqnarray}
If the inequality (\ref{approx19}) holds, then by the inequality $\frac{|x-y|^2}{\delta}+\delta\geq 2 |x-y|$, for every $\delta>0$, considering $(\o{x},\o{y},\o{t})\in \R^N\times (0,T)$, with $\o{x}\neq \o{y}$, 
and choosing $\delta=|\o{x}-\o{y}|$, by (\ref{approx19}) it follows,
\begin{eqnarray}\label{approx4}
Q(\o{x},\o{t})-Q(\o{y},\o{t})\leq 2M |\o{x}-\o{y}|+\var|\o{x}|^2+\frac{\gamma}{T-\o{t}}.
\end{eqnarray}
Then letting $\gamma,\;\var\rightarrow 0$, we have the assertion of Theorem \ref{approx1}.
Therefore it suffices to prove (\ref{approx19}). Suppose that (\ref{approx19}) were false. Then there would exist
$\delta_0,\;\gamma_0>0$, such that,
\begin{eqnarray}\label{approx4}
\sup_{\R^N\times[0,T)}H>0,\;\;\;\;\;\delta=\delta_0,\;\gamma=\gamma_0,
\end{eqnarray}
holds for a subsequence $\var=\var_n\rightarrow 0$. Since that $Q$ is bounded, we see $H<0$ for
sufficiently large $\|(x,y)\|$. Since that $2 M> Lip(v_0)$ we also see $H\leq 0$ at $t=0$.
Clearly $H\rightarrow -\infty$ at $t=T$; so (\ref{approx4}) now implies that $H$ takes its positive
maximum over $\R^N\times[0,T)$ at a point $(\w{x},\w{y},\w{t})$, with $\w{t}\in (0,T)$. First we study the 
behavior of the maximum point as $n\rightarrow \infty$. Since $H(\wx,\wy,\wt)>0$, by (\ref{approx4}), it follows from this obervation, that
\begin{eqnarray}\label{approx5}
\var|\wx|=\sqrt{\var^2 |\wx|^2}=\sqrt{\var}\sqrt{\var|\wx|^2}\leq \sqrt{\var (d-c)}=\sqrt{\var_n (d-c)}\rightarrow 0,\;n\rightarrow \infty.
\end{eqnarray}
Since $H$ attains its maximum over $\R^N\times[0,T)$ at $(\wx,\wy,\wt)$, if 
$w(x,y,t):=[Q(x,t)-\var|x|^2]-Q(y,t)$, we infer that,
\begin{eqnarray}\label{approx6}
(\partial_t K(\wx,\wy,\wt),\nabla K(\wx,\wy,\wt),\nabla^2 K(\wx,\wy,\wt))\in \P^{2,+}w(\wx,\wy,\wt).
\end{eqnarray}
Now we apply the usual Theorem of M.G. Crandall, and H. Ishii in \cite{2} , with $u_1(x,t):=Q(x,t)-\var |x|^2$, $u_2(y,t)=-Q(y,t)$; for
$\epsilon=\frac{\delta_0}{2 M}$, there are $(b_1,X_1)$, $(b_2,X_2)$ $\in\R\times\S^{N}$, such that,
\begin{eqnarray}\label{approx7}
(b_1, p+2\var \wx, X_1+2\var I_N)\in \o{\P}^{2,+}Q(\wx,\wt),\\
(-b_2, p, -X_2)\in \o{\P}^{2,-}Q(\wy,\wt),
\end{eqnarray}
where $p:=\frac{2 M}{\delta_0}(\wx-\wy)$. Moreover 
\begin{eqnarray}\label{approx9}
\left(\begin{array}{cc}
X_1 & 0\\
0 & X_2
\end{array}\right)\leq \frac{6 M}{\delta_0}\left(\begin{array}{cc}
I_N & -I_N\\
-I_N & I_N
\end{array}\right),\;\;\;b_1+b_2=\frac{\gamma_0}{(T-\wt)^2}.
\end{eqnarray}
To simplify notations we respectively denote with an over-bar, and an under-bar, the value of every fun\-ction which  
at $(\wx,\wt)$ and $(\wy,\wt)$, while for every function which depends only on the time
we omit that dependence.
By (\ref{approx8}), (\ref{approx7}), we can infer,
\begin{eqnarray}\label{approx10}
b_1 &\leq& \frac{1}{2}tr\big(\s\s^{\top}(X_1+2\var I_N)\big)+\la \o{\mu}, p+2\var \wx \ra
-e^{C\wt}\l_1(\o{Q} e^{C\wt})|\s^{\top}(p+2\var \wx)|^2\nonumber\\
&&-\l_2(\o{Q} e^{C\wt})
\la\s^{\top}(p+2\var \wx),\o{w}\ra-e^{-C\wt}\o{f}-C\o{Q},\nonumber\\ 
 b_2 &\leq& \frac{1}{2}tr\big(\s\s^{\top}X_2\big)-\la \u{\mu}, p \ra
+e^{C\wt}\l_1(\u{Q} e^{C\wt})|\s^{\top}p|^2\nonumber\\
&&+\l_2(\u{Q} e^{C\wt})
\la\s^{\top}p,\u{w}\ra+e^{-C\wt}\u{f}+C\u{Q}.\nonumber\\ \nonumber\\
\end{eqnarray}
Adding the inequalities (\ref{approx10}) and using (\ref{approx9}), yields,
\begin{eqnarray}\label{approx11}
\frac{\gamma_0}{T^2}&\leq& \frac{1}{2}tr\big(\s\s^{\top}(X_1+X_2+2\var I_N)\big)
+\Big[\la \o{\mu}-\u{\mu}, p \ra +2\la \o{\mu},\var \wx \ra\Big]\nonumber\\
&&+\Big[\N(\u{Q};p,\wx,\wt)-\N(\o{Q};p,\wx,\wt)\Big]
+\Big[-4e^{C\wt}\l_1(\o{Q} e^{C\wt})\la\s^{\top}p, \s^{\top}\var\wx \ra\nonumber\\
&&-4
e^{C\wt}\l_1(\o{Q} e^{C\wt})|\s^{\top}\var\wx|^2-2\l_2(\o{Q} e^{C\wt})
\la\s^{\top}\var \wx,\o{w}\ra\Big]\nonumber\\
&&+\Big[\l_2(\u{Q} e^{C\wt})\la\s^{\top}p,\u{w}-\o{w}\ra\Big]
+\Big[\u{f}-\o{f}\Big]e^{-C\wt}-C(\o{Q}-\u{Q}),
\end{eqnarray}
where, 
\begin{eqnarray}\label{approx12}
\N(Q;p,x,t):=e^{Ct}\l_1(Q e^{Ct})|\s^{\top}(t)p|^2+\l_2(Q e^{Ct})\la\s^{\top}(t)p,w(x,t)\ra.
\end{eqnarray}
Now, we estimate the single terms in the brackets $[\cdot]$ for the inequality (\ref{approx11}). 
Since $\l^{\p}_1>0$ over $[c,d]$, by (\ref{approx12}), we have,
\begin{eqnarray}\label{approx13}
\frac{ d\N}{dQ}&=& e^{2Ct}\l_1^{\p}(Q e^{Ct})|\s^{\top}(t)p|^2+e^{Ct}\l^{\p}_2(Q e^{Ct})
\la\s^{\top}(t)p,w(x,t)\ra\nonumber\\
&\geq& -\frac{|\l_2^{\p}(Q e^{Ct})|^2|w(x,t)|^2}{4\l_1^{\p}(Qe^{Ct})}
\geq -\frac{\|\l^{\p}_2\|^2_{\infty}\|w\|^2_{\infty}}{4\min_{[c,d]}\l_1^{\p}},\;\forall\;Q,\;p,\;x,t.
\end{eqnarray}
So by $\o{Q}-\u{Q}>0$
\begin{eqnarray}\label{approx13}
\N(\u{Q};p,\wx,\wt)-\N(\o{Q};p,\wx,\wt)\leq \frac{\|\l^{\p}_2\|^2_{\infty}\|w\|^2_{\infty}}{4\min_{[c,d]}\l_1^{\p}}(\o{Q}-\u{Q}).
\end{eqnarray}
>From inequality (\ref{approx9}), we have,
\begin{eqnarray}\label{approx14}
\frac{1}{2}tr\big(\s\s^{\top}(X_1+X_2+2\var I_N)\big)\leq \var tr(\s\s^{\top}).
\end{eqnarray}
Conditions (\ref{approx4}) and (\ref{approx5}) yield,
\begin{eqnarray}\label{approx15}
\la \o{\mu}-\u{\mu}, p \ra +2\la \o{\mu},\var \wx \ra &\leq& 
Lip(\mu)\frac{2M}{\delta_0}|\wx-\wy|^2+2\|\mu\|_{\infty}\sqrt{\var (d-c)}\nonumber\\
&\leq& 2 Lip(\mu)(\o{Q}-\u{Q})+2\|\mu\|_{\infty}\sqrt{\var (d-c)},
\end{eqnarray}
\begin{eqnarray}\label{approx16}
\l_2(\u{Q} e^{C\wt})\la\s^{\top}p,\u{w}-\o{w}\ra\leq 2\|\l_2\|_{\infty}\|\s^{\top}\|_{\infty}Lip(w)(\o{Q}-\u{Q}).
\end{eqnarray}
\begin{eqnarray}\label{approx16*}
\Big[\u{f}-\o{f}\Big]e^{-C\wt}&\leq& Lip(f)\Big[|\wx-\wy|+(\o{Q}-\u{Q})e^{C \wt}\Big]e^{-C\wt}\nonumber\\
&\leq& Lip(f)\Big[\frac{1}{2}\Big(\frac{|\wx-\wy|^2}{\delta_0}
+\delta_0\Big)+(\o{Q}-\u{Q})e^{C \wt}\Big]e^{-C\wt}\nonumber\\
&\leq& Lip(f)\Big[\frac{\o{Q}-\u{Q}}{2 M}+(\o{Q}-\u{Q})e^{C \wt}\Big]e^{-C\wt}\nonumber\\
& = & Lip(f)(\o{Q}-\u{Q})\Big(\frac{e^{-C \wt}}{2 M}+1\Big)\nonumber\\
&\leq & Lip(f)\Big(\frac{1}{2 M}+1\Big)(\o{Q}-\u{Q}). 
\end{eqnarray}
\end{Dimo}
Being $\delta_0$ fixed, $p$ is bounded for $\var\rightarrow 0$; then
by (\ref{approx5}), the term $[-4e^{C\wt}\times\l_1(\o{Q} e^{C\wt})\la\s^{\top}p, \s^{\top}\var\wx \ra-4
e^{C\wt}\l_1(\o{Q} e^{C\wt})|\s^{\top}\var\wx|^2-2\l_2(\o{Q} e^{C\wt}) 
\la\s^{\top}\var \wx,\o{w}\ra\Big]$ in (\ref{approx11}), is of order $\sqrt{\var}$, for $\var\rightarrow 0$. 
Using that observation and introducing the estimates (\ref{approx13})-(\ref{approx16*}), in the inequality (\ref{approx11}), yields
\begin{eqnarray}\label{approx17}
\frac{\gamma_0}{T^2}&\leq & \var tr(\s\s^{\top})
+2 Lip(\mu)(\o{Q}-\u{Q})+2\|\mu\|_{\infty}\sqrt{\var (d-c)}\nonumber\\
&&+\frac{\|\l^{\p}_2\|^2_{\infty}\|w\|^2_{\infty}}{4\min_{[c,d]}\l_1^{\p}}(\o{Q}-\u{Q})
+O(\sqrt{\var})
+2\|\l_2\|_{\infty}\|\s^{\top}\|_{\infty}Lip(w)(\o{Q}-\u{Q})\nonumber\\
&&+Lip(f)\Big(\frac{1}{2 M}+1\Big)(\o{Q}-\u{Q})-C(\o{Q}-\u{Q})
\end{eqnarray}
Letting $n\rightarrow \infty$ in (\ref{approx17}), we see
\begin{eqnarray}\label{approx18}
\frac{\gamma_0}{T^2}&\leq &
\Big[2 Lip(\mu)
+\frac{\|\l^{\p}_2\|^2_{\infty}\|w\|^2_{\infty}}{4\min_{[c,d]}\l_1^{\p}}
+2\|\l_2\|_{\infty}\|\s^{\top}\|_{\infty}Lip(w)+\nonumber\\
&& Lip(f)\Big(\frac{1}{2 M}+1\Big)-C\Big](\o{Q}-\u{Q})=0.
\end{eqnarray}
In the last passage we have used the definition (\ref{Constant}). The inequality (\ref{approx18}) contradicts $\gamma_0>0$. We 
thus prove (\ref{approx19}).\cvd


\end{document}